\documentclass{book}    
\usepackage{piers}  

\usepackage{algorithm}
\usepackage{algpseudocode}
\usepackage{float}
\usepackage{caption}
\usepackage{subcaption}
\usepackage{graphicx}
\DeclareMathOperator\supp{supp}
\DeclareMathOperator\rect{rect}
\DeclareMathOperator\DOF{DOF}
\usepackage[colorinlistoftodos,textwidth=4cm,shadow]{todonotes}

\newcommand{\mi}{\mathrm{i}}

\begin{document}

\title{Vico-Greengard-Ferrando quadratures in the tensor solver for integral equations}
\maketitle

\begin{authors}

{\bf V. Khrulkov}$^{1}$, {\bf M. Rakhuba}$^{1}$, {\bf I. Oseledets}$^{1,2}$\\
\medskip
$^{1}$Skolkovo Institute of Science and Technology, Russia\\
$^{2}$Institute of Numerical Mathematics, Russia

\end{authors}

\begin{paper}

\begin{piersabstract}
Convolution with Green's function of a differential operator appears in a lot of applications e.g. Lippmann-Schwinger integral equation. Algorithms for computing such are usually non-trivial and require non-uniform mesh. However, recently Vico, Greengard and Ferrando developed method for computing convolution with smooth functions with compact support with spectral accuracy, requiring nothing more than Fast Fourier Transform (FFT). Their approach is very suitable for the low-rank tensor implementation which we develop using Quantized Tensor Train (QTT) decomposition. 
\end{piersabstract}

\psection{Introduction}
In this paper we propose an algorithm for computing an approximation to a convolution
$$ h(\textbf{r}') = \int_{\mathbb{R}^3} g(\textbf{r}-\textbf{r}') \rho(\textbf{r}) d \textbf{r}, $$
where $G$ is a continuous Green's function of some PDE. It has been shown recently in \cite{vico2016fast} that discrete approximation to $h$ yieding spectral accuracy for compactly supported smooth functions $\rho$ can be constructed using relatively simple idea. Idea is based on the fact that continuous Green's function can be replaced by the truncated one, Fourier transform of which is in $\mathcal{C}^{\infty}(\mathbb R ^{3})$. Applying convolution theorem one gets answer as a Fourier type integral of a rapidly decreasing function from $\mathcal{C}^{\infty}(\mathbb R ^{3})$ which is then discretized on finer grid and evaluated using Discrete Fourier Transform (DFT). This scheme achieves spectral order of accuracy due to the superalgebraic approximation of the continuous Fourier transform by DFT and has been shown to also give second order with small constant for continuos functions. We implemented this scheme using data compression via Quantized Tensor Train (QTT). Approach is straightforward -- we compute discrete truncated Green's function using formulas from \cite{vico2016fast} and convert it to QTT-format. Then we use standard algorithms in QTT which allow for logarithmic storage with respect to the grid size and accuracy of approximation. Moreover, discrete convolution and solution of integral equations are then performed also with logarithmic complexity by the algorithms described in \cite{khkaz-conv-2013},\cite{osel-2d2d-2010}. 
Paper is organized as follows:
\begin{itemize}
    \item We briefly describe the algorithm developed in Section~\ref{sec1}
    \item We present necessary definitions and algorithms from TT theory in Section~\ref{sec-tt}
    \item We describe TT implementation of the algorithm in Section~\ref{sec-qtt}
    \item We present numerical results in Section~\ref{sec:num}
\end{itemize}
The main reference for Vico-Greengard-Ferrando quadrature is the original paper \cite{vico2016fast}, in the presentation of the algorithm we follow \cite{klinteberg2016fast}.
\psection{Algorithm} \label{sec1}
We wil present the scheme described in \cite{vico2016fast} in a way which is most suitable for transitioning to TT-format. 
Let $\rho(\textbf r)$ be a smooth function such that $\supp \rho \subset D = [0,L_1] \times [0, L_2] \times [0, L_3] $ and we are interested in computing
\begin{equation}
h(\textbf r) = \int_{\mathbb{R}^3} g(\textbf r-\textbf r') \rho(\textbf r') d \textbf r',
\label{eq1}
\end{equation}
where $g(\textbf r)$ is a Green's function of some differential operator. Let us assume that $g(\textbf r)$ depends only on $r = |\textbf r|$. If we seek restriction of the solution $h(\textbf r)$ to $D$ then $|r-r'|$ in $\eqref{eq1}$ doesn't exceed $L = \sqrt {L_1^2 + L_2 ^2 + L_3 ^2} $. Thus if we replace $g(\textbf r)$ by $g_L(\textbf r) = g(\textbf r) \rect{\frac{r}{2L}}$,
$$\int_{\mathbb{R}^3} g(\textbf r-\textbf r') \rho(\textbf r') d\textbf r' = \int_{\mathbb{R}^3} g_L(\textbf r-\textbf r') \rho(\textbf r') d\textbf r'.$$
Advantage of using $g_L(\textbf r)$ is that since it has compact support it's Fourier transform $G_L(\textbf s)$ is in $C^{\infty} (\mathbb{R}^3)$, and is straightforward to compute for many differential operators. We mainly focus on Helmholtz differential operator $\nabla ^2 + k^2$ for which the following formula holds:
\begin{gather}
	g(\textbf r) = \frac{e^{\mi kr}}{4 \pi r}, \notag \\
	G_L(\textbf s) = \int_{\mathbb{R}^3} g(\textbf r) \rect{\left(\frac{r}{2L}\right)} e^{\mi \textbf s \cdot \textbf r} d \textbf r \notag \\
	= \frac{-1 + e^{\mi Lk} (\cos Ls - \mi \frac{k}{s} \sin Ls)}{(k-s)(k+s)}.
\end{gather}
It is easy to check that $G_L(\textbf s)$ is indeed smooth and nonsingular.
If we denote Fourier transform of $\rho(\textbf r)$ by $\widehat{\rho}(\textbf s)$ we obtain the final formula:
\begin{equation}
h( \textbf r ) = \Big(\frac{1}{2 \pi} \Big)^3\int_{\mathbb{R}^3} e^{\mi \textbf s \cdot \textbf r} \widehat{\rho}( \textbf s) G_L(\textbf s)d \textbf s.
\label{eq2}
\end{equation}
This integral is then discretized using trapezoidal rule and computed using DFT. However to cancel slightly oscillatory behavior of $\widehat{G_L}(s)$ zero padding by a factor of at least $3$ is required (we will use factor $4$ to keep grid size being a power of $2$ - for the analysis see \cite{klinteberg2016fast}). 
Suppose that domain $D$ is discretized using uniform grid with  $N_i$ nodes in corresponding dimensions and function $\rho$ is sampled on this grid yielding an array $\rho _{ijk}$ i.e.
$$\rho_{ijk} := \rho(i h_1, j h_2,k h_3), \quad h_i = \frac{L_i}{N_i}.$$
Algorithm then is summarized as following
\begin{algorithm}
\caption{Basic Vico-Greengard-Ferrando quadrature}
\label{basic}
\begin{algorithmic}[1]
\State Zero pad $\rho_{ijk}$ by a factor of $4$ and compute 3d FFT 
defining $\widehat{\rho}({\textbf s})$ for $\textbf{s} =\frac{\pi}{2} (\frac{s_1}{L_1},\frac{s_2}{L_2},\frac{s_3}{L_3})$,
$s_i \in \lbrace -2N_i, \hdots, 2N_i-1 \rbrace$.

\State Evaluate $G_L(\textbf{s})$ for $L = \sqrt{L_1^2 + L_2^2 + L_3^2}$ for $\textbf{s}$ defined above, and multiply elementwise by $\widehat{\rho}({\textbf s})$. 
\State Perform 3d IFFT on the array defined above and truncate the result keeping first $N_1 \times N_2 \times N_3$ entries, obtaining approximation to $h(\textbf r)$ on the grid.
\
\end{algorithmic}
\end{algorithm}
\\

\psubsection{Convolution form of the algorithm}
Let us write Algorithm \ref{basic} more explicitly. 
All steps together can be represented as follows
\begin{equation}\label{ifftbig}
h_{ijk} = \frac{1}{4N_1}  \frac{1}{4N_2}  \frac{1}{4N_3} \sum_{s_1 = - 2 N_1}^{2 N_1 -1} \sum_{s_2 = - 2 N_2}^{2 N_2 -1} \sum_{s_3 = - 2 N_3}^{2 N_3 -1} G_L(\textbf s) \widehat{\rho}({\textbf s}) e^{2 \pi \mi \frac{s_1 i}{4 N_1}} e^{2\pi \mi \frac{s_2 j}{4 N_2}} e^{2\pi \mi \frac{s_3 k}{4 N_3}} ,
\end{equation}
where 
\begin{equation}\label{fftbig}
\widehat{\rho} ({\textbf s}) = \sum_{i' = 0}^{4N_1 - 1}\sum_{j' = 0}^{4N_2 - 1} \sum_{k' = 0}^{4N_3 - 1} 
\rho_{i'j'k'} e^{-2\pi \mi \frac{s_1 i'}{4N_1}}e^{-2\pi \mi \frac{s_2 j'}{4N_2}}  e^{-2\pi \mi \frac{s_3 k'}{4N_3}}.
\end{equation}  
By plugging \eqref{fftbig} into \eqref{ifftbig} and by changing the order of summation it is easy to see that
\begin{equation}\label{convolution}
h_{ijk} = \sum_{i' = 0}^{4N_1 - 1}\sum_{j' = 0}^{4N_2 - 1} \sum_{k' = 0}^{4N_3 - 1} G^M_{i-i',j-j',k-k'} \rho_{i'j'k'},
\end{equation}
where 
\begin{equation} \label{eq:GM}
G^M_{i-i',j-j',k-k'} = \frac{1}{4N_1}  \frac{1}{4N_2}  \frac{1}{4N_3} \sum_{s_1 = - 2 N_1}^{2 N_1 -1} \sum_{s_2 = - 2 N_2}^{2 N_2 -1} \sum_{s_3 = - 2 N_3}^{2 N_3 -1} G_L(\textbf s) e^{2 \pi \mi \frac{s_1 (i-i')}{4 N_1}} e^{2\pi \mi \frac{s_2 (j-j')}{4 N_2}} e^{2\pi \mi \frac{s_3 (k-k')}{4 N_3}} .
\end{equation}
Moreover, since $\rho_{i'j'k'}$ is $0$  for $ i' \geq N_1, j' \geq N_2, k' \geq N_3$ (see step $1$ of Algorithm~\ref{basic}) and we truncate the result, formula \eqref{convolution} simplifies and finally:
\begin{equation}\label{convolution-final}
h_{ijk} = \sum_{i' = 0}^{N_1 - 1}\sum_{j' = 0}^{N_2 - 1} \sum_{k' = 0}^{N_3 - 1} G^M_{i-i',j-j',k-k'} \rho_{i'j'k'}.
\end{equation}
We see that \eqref{convolution-final} takes the form of a discrete aperiodic convolution with discrete Green's function $G^M$ (which we will call \textit{mollified Green's function}). One can notice that to fully determine $G^M$ it is sufficient to run Algorithm \ref{basic} once for a special right hand side $\rho_{ijk} = \delta_{i0} \delta_{j0}\delta_{k0}$.
Formula \eqref{convolution-final} plays essential role in the further analysis. Multiplication by multilevel Toeplitz matrix generated by $G^M$ can be performed with logarithmic complexity in QTT format as described in \cite{khkaz-conv-2013}, and we discuss neccessary definitions and algorithms in the next section.

\psection{Low-rank tensor approach} \label{sec2}

\psubsection{TT and QTT formats} \label{sec-tt}
To understand the QTT format let us start with describing the TT-format, which is a nonlinear low-paramentric representation of multidimensional arrays, called \emph{tensors}.
Tensor $\mathcal{X}\in\mathbb{C}^{n_1\times n_2 \times \dots \times n_d}$ is said to be in the TT-format if it represents as
\begin{equation}\label{tt-dec}
\mathcal{X}_{i_1 i_2\dots i_d} = X^{(1)}(i_1) X^{(2)}(i_2) \dots X^{(d)}(i_d),
\end{equation}
where $X^{(k)}(i_k)\in \mathbb{C}^{r_{k-1}\times r_k}$, $r_0=r_d=1$, $i_k=1,\dots,n_k$.
Matrices $X^{(k)}$ are called \emph{TT-cores} and $r_k$ are called \emph{TT-ranks}.
Notice that if $r = \max_k r_k$ is small, then there is a significant compression to store $\mathcal{X}$.
Indeed, initial tensor requires storing $n^d$ parameters, while to store its TT-representation only $\mathcal{O}(dnr^2)$ parameters are needed.

In fact, one could use TT representation to store and to work with arising in Algorithm~\ref{basic} 3-dimensional arrays. However we will use a more sophisticated approach called QTT format, which allows for additional storage reduction compared to TT.
QTT format is the following modification of the TT format.
First we assume that $d=3$, $n_i=2^{l_i}$, $i=1,2,3$.
Then each ``physical'' index $i$, $j$, $k$ is represented in the binary format, i.e. 
$$
	i = i_1 + 2^1\, i_2 + \dots + 2^{l_i -1} i_d, \quad i_m = 0,1, \quad m=1,\dots,d
$$
and we have initial tensor $\rho_{ijk}$ encoded as a $(l_1+l_2 + l_3)$-dimensional array $\tilde\rho$:
$$
	\rho_{ijk} \equiv \tilde\rho_{i_1\dots i_{l_1} j_1\dots j_{l_2} k_1\dots k_{l_3}}.
$$
TT decomposition of $\tilde\rho$ is called the QTT decomposition.
The storage of the QTT decomposition is $\mathcal{O}(r^2(l_1+l_2 + l_3)) = \mathcal{O}(r^2 \log n)$, so if $r=\max_i r_i$ is bounded, the total storage scales logarithmically.
In practice tensors of exact low-rank rarely occur.
Typically one fixes accuracy $\epsilon$ and tries to find best approximation with this accuracy. It has been shown that in some applications ranks grow as $r=\mathcal{O}(\log^\alpha \epsilon^{-1})$, $\alpha>1$ \cite{ks-2dqtt-2015,ksro-multiscale-2016}.

\psubsection{Translation of the algorithm to the QTT format} \label{sec-qtt}


\paragraph{Computation of kernel $G^M$.} \label{kern-comp}
To use \eqref{convolution-final} we first need to find $G^M$ \eqref{eq:GM}. For this purpose we run Algorithm~\ref{basic} for $\rho_{ijk} = \delta_{i0}\delta_{j0}\delta_{k0}$.
Precomputations are done in the full format, in other words we form the whole dense tensor $G^M$ and utilize TT-SVD algorithm \cite{osel-tt-2011} to find its QTT representation.
TT-SVD algorithm is based on the computation of SVD decompositions of tensors reshaped into full 2D matrices and therefore is quite expensive.
In principle one could use DFT in the QTT format \cite{dks-ttfft-2012}  to avoid forming full tensors.
Unfortunately, we found that intermediate tensors arising in Algorithm~\ref{basic} are of large rank. 
We will address this problem in our future work.

\paragraph{Computation of $\rho$.}
Tensor $\rho$ can be already given in the QTT representation.
This can happen, e.g. if we are running a certain iterative process involving computation of convolution \eqref{eq1} and all operations in this process are done within the QTT format.
Otherwise, $\rho_{ijk}$ can be approximated with logarithmic complexity by using the cross approximation method \cite{ot-ttcross-2010}, which adaptively samples elements of a tensor.
In this case we just need $\rho_{ijk}$ be given as a function which returns value by given 3 indices $i,j,k$.

\paragraph{Computation of convolution $G^M * \rho$.}
Next goal is to find convolution of tensors $G^M$ and $\rho$ \eqref{convolution-final}.
The convolution can be considered as a multiplication of multilevel Toeplitz matrix generated by $G^M$ and vector $\rho$.
Matrices can also be represented in the TT and by analogy in the QTT format. 
The definition is similar to that of TT-tensor: given matrix (operator) $\mathcal{A}_{i_1\dots i_d j_1\dots j_d}$, which acts on vector $\mathcal{X}_{j_1\dots j_d}$ such that
$$
	\mathcal{Y}_{i_1\dots i_d} = \sum_{j_1\dots j_d} \mathcal{A}_{i_1\dots i_d j_1\dots j_d}\mathcal{X}_{j_1\dots j_d}.
$$
its TT-decomposition is defined as
\[
	\mathcal{A}_{i_1\dots i_d j_1\dots j_d} = A^{(1)}(i_1,j_1) A^{(2)}(i_2,j_2) \dots A^{(d)}(i_d,j_d),
\]
where $A^{(k)}(i_k,j_k)\in \mathbb{C}^{R_{k-1}\times R_k}$, $R_0=R_d=1$, $i_k=1,\dots,n_k$.
QTT decomposition of 3D operator $G^M_{i-i',j-j',k-k'}$ is defined by analogy with the QTT decomposition of a tensor -- we quantize indices $i,j,k$ and $i',j',k'$, group them pairwise and then compute TT decomposition:
\[
\begin{split}
	G^M_{i-i',j-j',k-k'} = & G^{(1)}(i_1,i_1')\dots  G^{(l_1+2)}(i_{l_1+2},i_{l_1+2}') \\
	& G^{(l_1+3)}(j_1,j_1')\dots  G^{(l_1+l_2+4)}(j_{l_2+2},j_{l_2+2}') \\
	& G^{(l_1+l_2 + 5)}(k_1,k_1')\dots  G^{(l_1+l_2+l_3+6)}(k_{l_3+2},i_{l_3+2}').
\end{split}
\]
We use approach from \cite{khkaz-conv-2013} and analytically construct QTT representation of the induced multilevel Toeplitz matrix $G^M_{i-i',j-j',k-k'}$ given QTT representation of $G^M_{i,j,k}$.
Then matrix-vector product \eqref{convolution-final} can be done in different ways.
We used optimization procedure AMEn (alternating minimal energy method) \cite{ds-amr1-2013,ds-amr2-2013} which allows for rank adaptation compared to standard ALS (alternating least squares) \cite{holtz-ALS-DMRG-2012} optimization which works with the representation of a given size.

\psection{Numerical experiments}\label{sec:num}
\paragraph{Approximating $G^M$ using QTT.}
Firstly we show that using QTT representation greatly reduces number of degrees of freedom (DOF) of $G^M$.
Suppose that tensor $\mathcal{X}$ is given in the QTT format with ranks $\{r_1,\dots,r_d\}$. Then it is easy to count total number of DOF of $\mathcal{X}$:
$$\DOF(\mathcal{X}) = 2(r_1 + r_d) + \sum_{i=2}^{d} 2 r_{i-1}r_i.$$
By applying this formula to $G^M$ computed as described in \ref{kern-comp} for various values of $k$ (while keeping $L=1$) we obtained the following results (see Figure~\ref{fig:dof}). This shows advantages of using QTT.
\begin{figure}[h]
	\centering
	\begin{subfigure}{.45\textwidth}
		\centering
		     	\includegraphics[width=\linewidth]{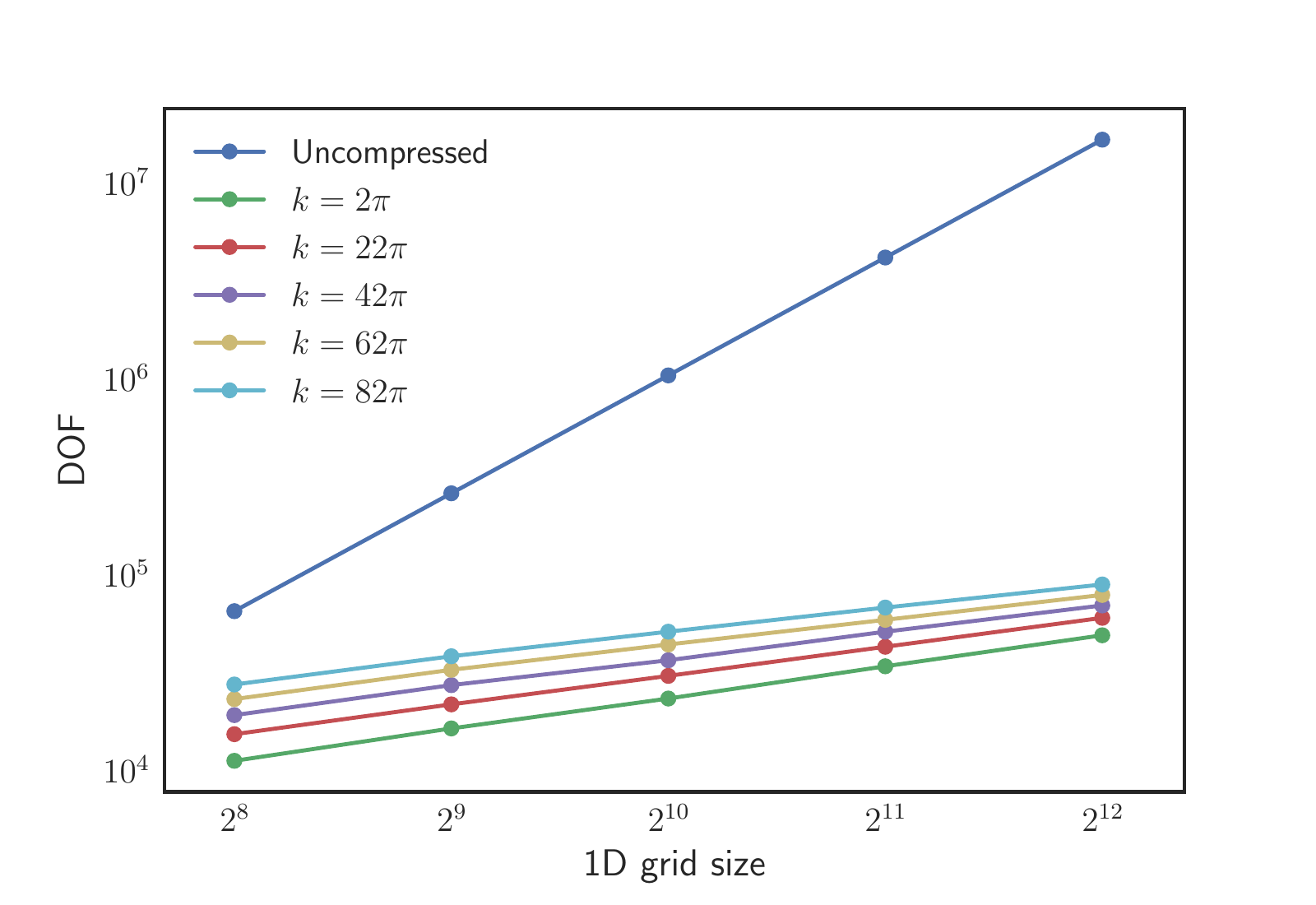}
		        \caption{ }
		     	\label{fig:dof2d}
	\end{subfigure}
		\begin{subfigure}{.45\textwidth} 
		\centering
		     	\includegraphics[width=\linewidth]{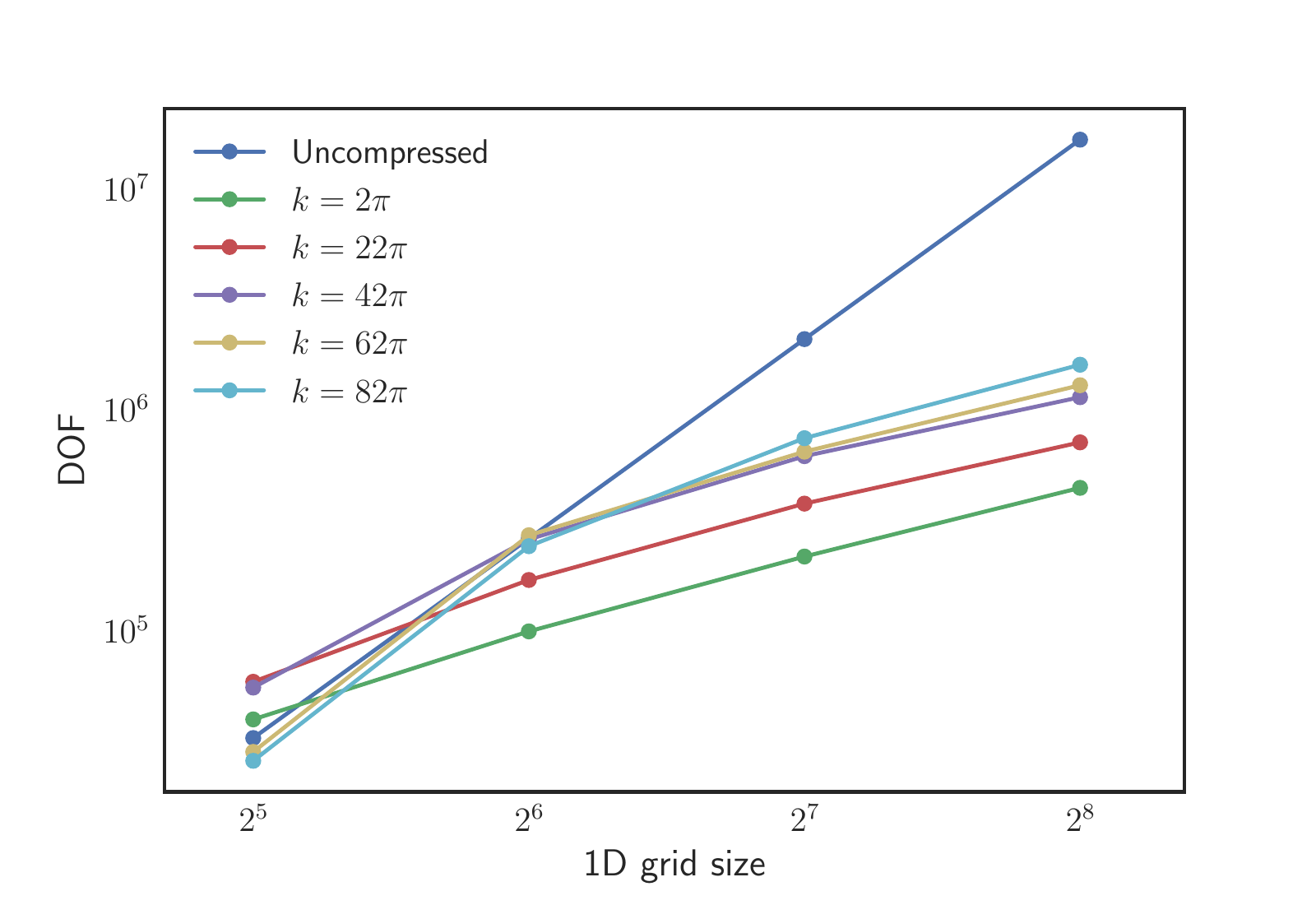}
		     	\caption{ }
		     	\label{fig:dof3d}
	\end{subfigure}%
	
	\caption{DOF of $G^M$ in $2D$ (\ref{fig:dof2d}) and $3D$ (\ref{fig:dof3d}) for $\epsilon = 10^{-7}$ and $L=1$}
	\label{fig:dof}
\end{figure}
\paragraph{Solving scattering problems.}
To further test our approach we solve the Lippmann-Schwinger equation which is used for solving scattering problems:
\begin{equation}\label{eq:LS}
\sigma (\textbf{r}) + k^2 q(\textbf{r}) \int_{\mathbb{R}^3} \frac{e^{\mi k | \textbf{r} - \textbf{r}'|}}{| \textbf{r} - {\bf r}'|} \sigma(\textbf{r}') d\textbf{r}' = -k^2 q(\textbf{r}) \phi^\mathrm{inc},
\end{equation}
and then we find
$$\phi^\mathrm{scat} = \int_{\mathbb{R}^3} \frac{e^{ik | \textbf{r} - \textbf{r}'|}}{| \textbf{r} - \textbf{r}'|} \sigma(\textbf{r}') d\textbf{r}'.$$
We used rounding by $\epsilon = 10^{-7}$ in our computations.
Firstly we fixed $k=1$ and $L=32 \pi$ and took $q(\textbf{r})$ to be a $3D$ gaussian:
\begin{equation}\label{eq:gaussian}
q(\textbf{r}) = e^{\frac{-|\textbf{r}-\textbf{r}'|^2}{2 a^2}},
\end{equation}
with $\textbf{r}' = \left(\frac{L}{2},\frac{L}{2},\frac{L}{2}\right),$ and $a = \frac{L}{10},$
and $$\phi^\mathrm{inc}(x,y,z) = e^{\mi x}.$$

To solve arising systems in the QTT format we used AMEn \cite{ds-amr1-2013,ds-amr2-2013} , which allows for rank adaptation.
\begin{figure}[h!]
\centering
	\begin{subfigure}{.45\textwidth}
		\centering
		     	\includegraphics[width=\linewidth]{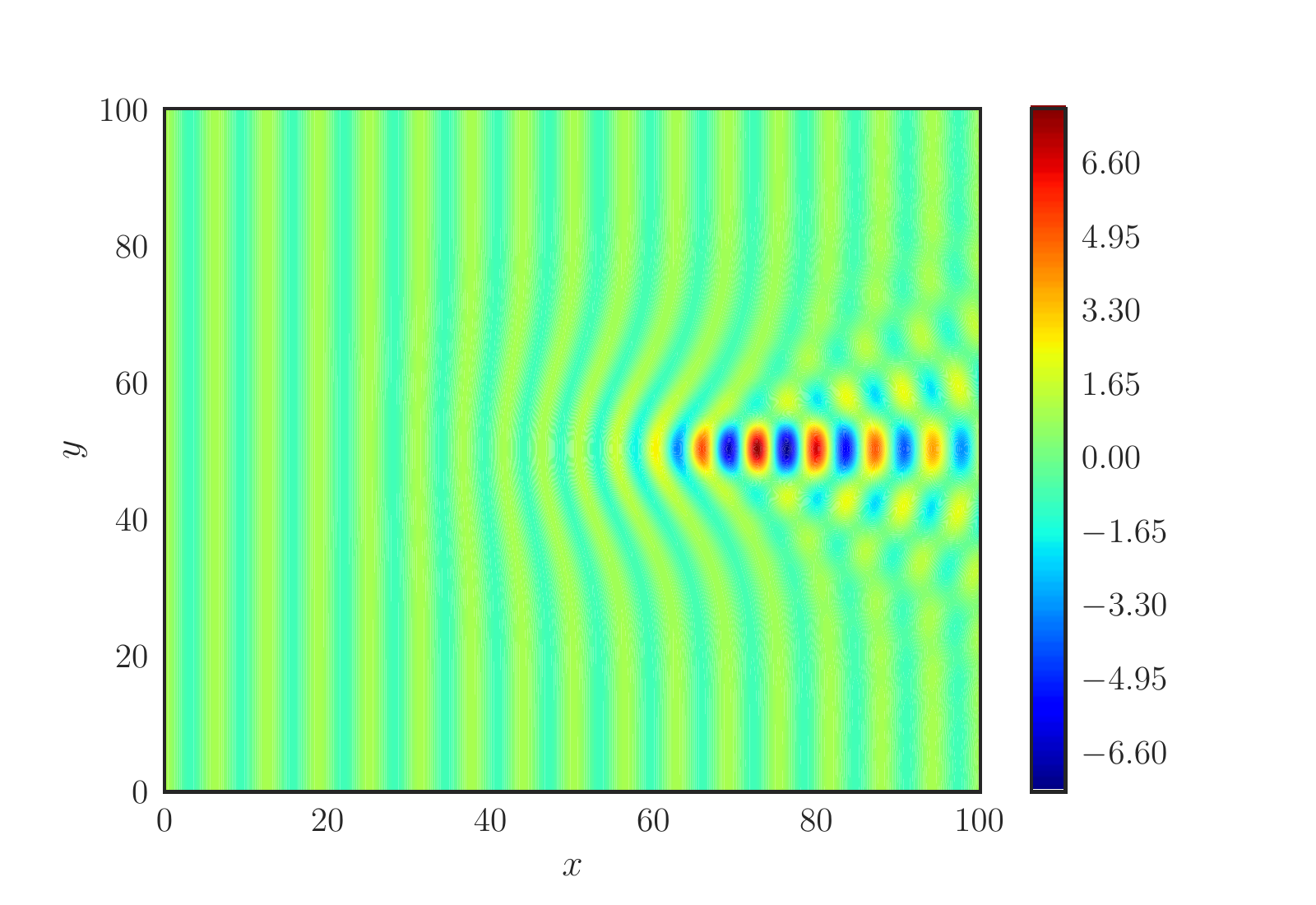}
	     \caption{Section $z = 16 \pi.$}
	     	
	\end{subfigure}
	\begin{subfigure}{.45\textwidth} 
		\centering
		     	\includegraphics[width=\linewidth]{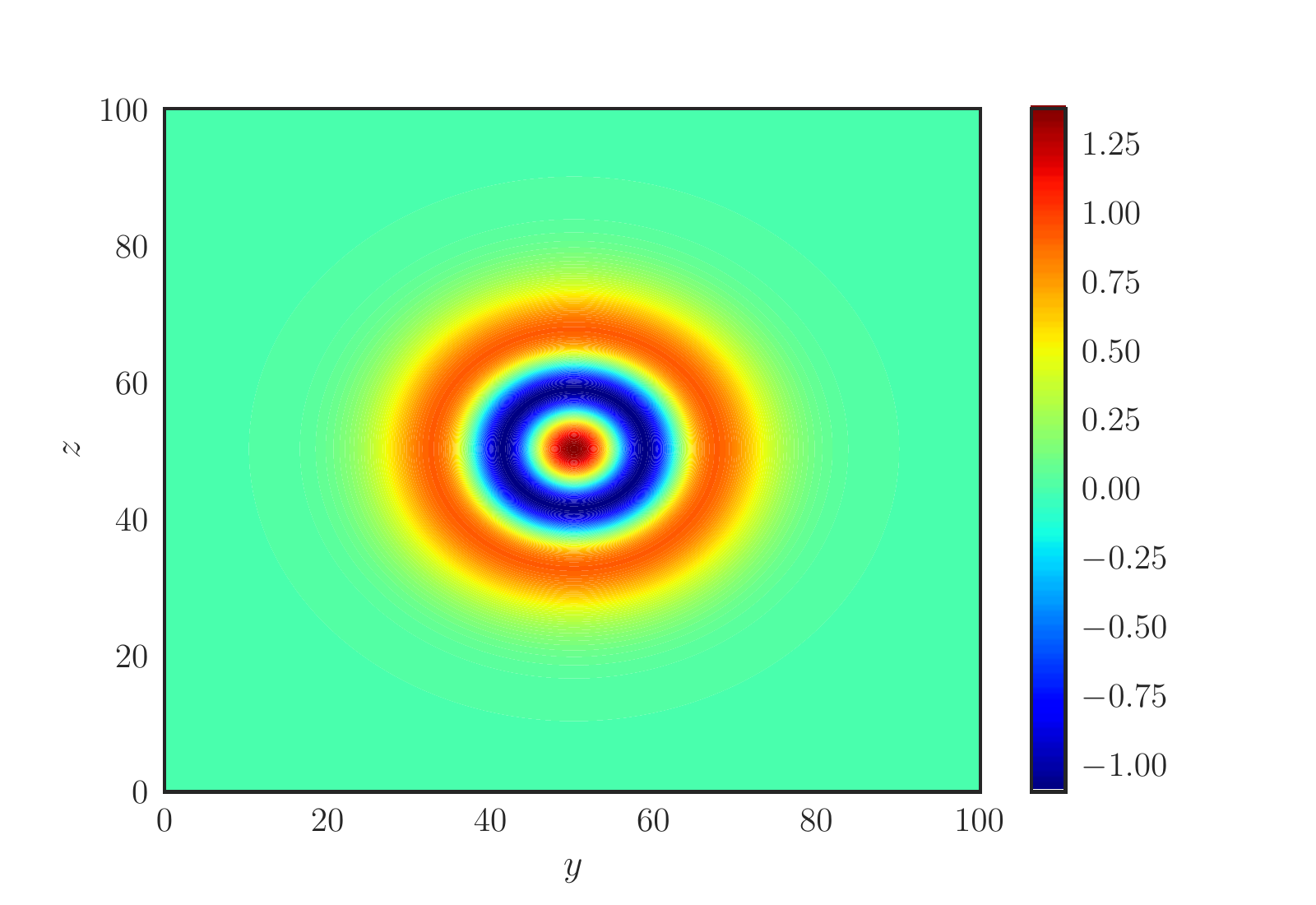}
		    \caption{Section $x = 17.5 \pi.$}
		     	
	\end{subfigure}%
	\caption{Scattering on the Gaussian with standard deviation $3.2\pi$ in $D = [0,32\pi] \times [0,32\pi] \times [0,32\pi] $ for $k=1$. Effective rank of the solution is equal to $39.$}
	\label{fig:gauss}
\end{figure}
Taking the solution computed on a grid with size $256^3$ as a reference we measured relative error of the solutions computed on smaller grids for various rounding parameters $\epsilon$. Results are given in the Table~\ref{table:err}. We see that error is roughly equal to the $\epsilon$ even for modest number of grid nodes per wavelength. 
\begin{table}
\caption{Relative error (err) and effective rank (erank) for different grid sizes and rounding errors $\epsilon$. Results are presented for two types of function $q$: Gaussian \eqref{eq:gaussian} and smoothed cube \eqref{eq:smoothed-cube}.}
\centering
\begin{tabular}[t]{lccccccc}
	\multicolumn{1}{c}{} & \multicolumn{1}{c}{3D grid size} & \multicolumn{2}{c}{$32^3$} & \multicolumn{2}{c}{$64^3$} & \multicolumn{2}{c}{$128^3$} \\
	\hline
	&  & err & erank &  err & erank & err & erank \\
	\hline
 Gaussian	& $\epsilon = 10^{-3}$ & 8e-1 &14 &8e-2 &15 &7e-2 &14 \\
 				& $\epsilon = 10^{-5}$ & 8e-1 &24 &4e-4 &25  &4e-4 &26 \\
 				& $\epsilon = 10^{-7}$ & 8e-1 & 30 & 5e-6 & 34 & 4e-6 & 38\\
 	\hline
  Smoothed cube	& $\epsilon = 10^{-3}$ &8e-1&22 &8e-2 &35 &6e-2 &34 \\
 				& $\epsilon = 10^{-5}$ &8e-1 &29 &5e-4 &57 &5e-4 & 59 \\
 				& $\epsilon = 10^{-7}$ &  8e-1 & 29 & 2e-5 & 73 & 6e-6 & 80\\
 	\hline
\end{tabular}
\label{table:err}
\end{table}

As a next experiment we performed the same computations for $q(r)$ representing smoothed cube:
\begin{equation}\label{eq:smoothed-cube}
q(\textbf{r}) = e^{-0.5 \left(\frac{| \textbf r-{\textbf r}'|}{a}\right)^8},
\end{equation} for 
$\textbf{r}' = \left(\frac{L}{2},\frac{L}{2},\frac{L}{2}\right)$ and $a = \frac{L}{4}$.
For the results see Figure~\ref{fig:cube}.
\begin{figure}[h!]
	\centering
	\begin{subfigure}{.45\textwidth}
		\centering
		     	\includegraphics[width=\linewidth]{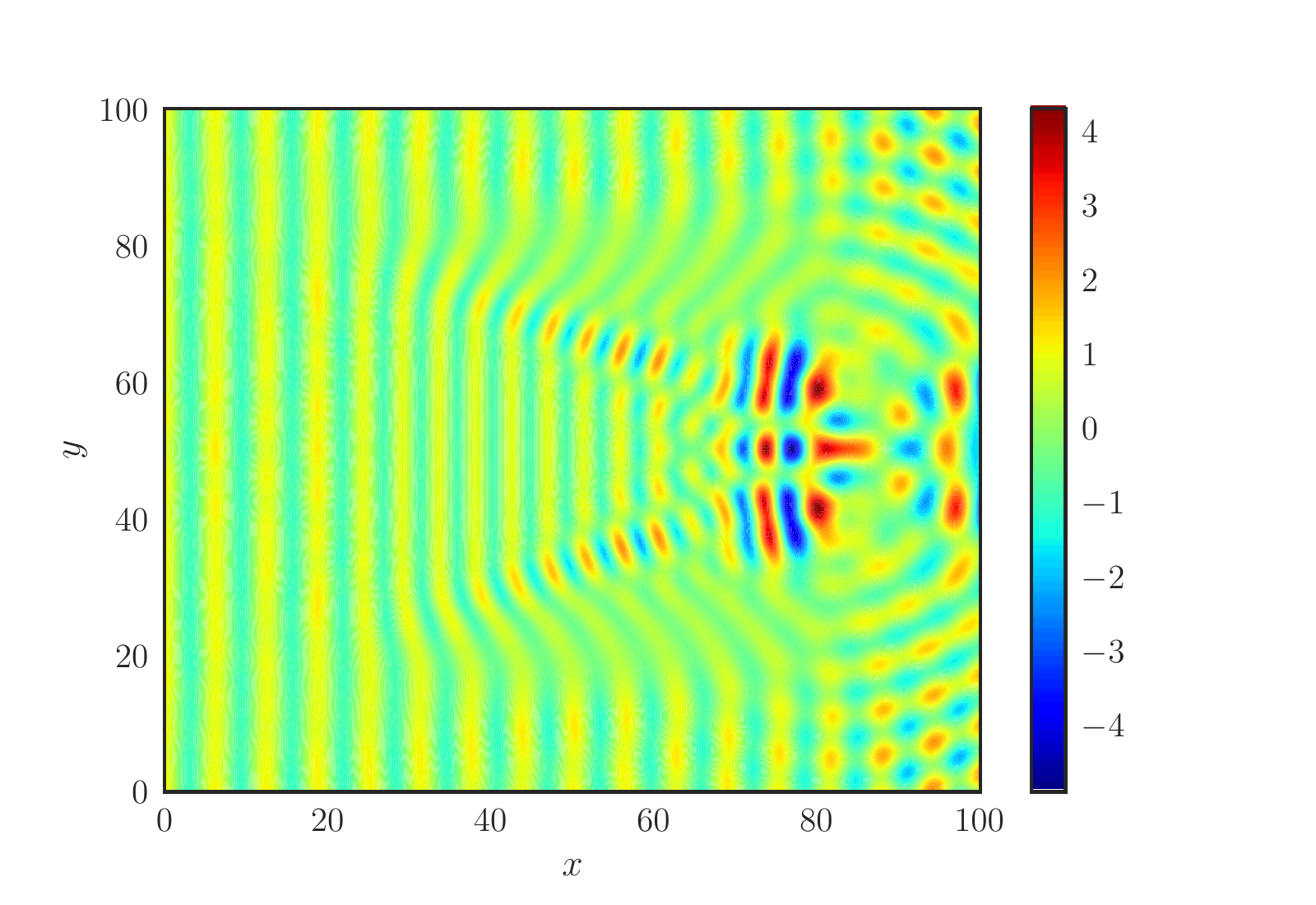}
		  	     \caption{Section $z = 16 \pi.$}
	\end{subfigure}
		\begin{subfigure}{.45\textwidth} 
		\centering
		     	\includegraphics[width=\linewidth]{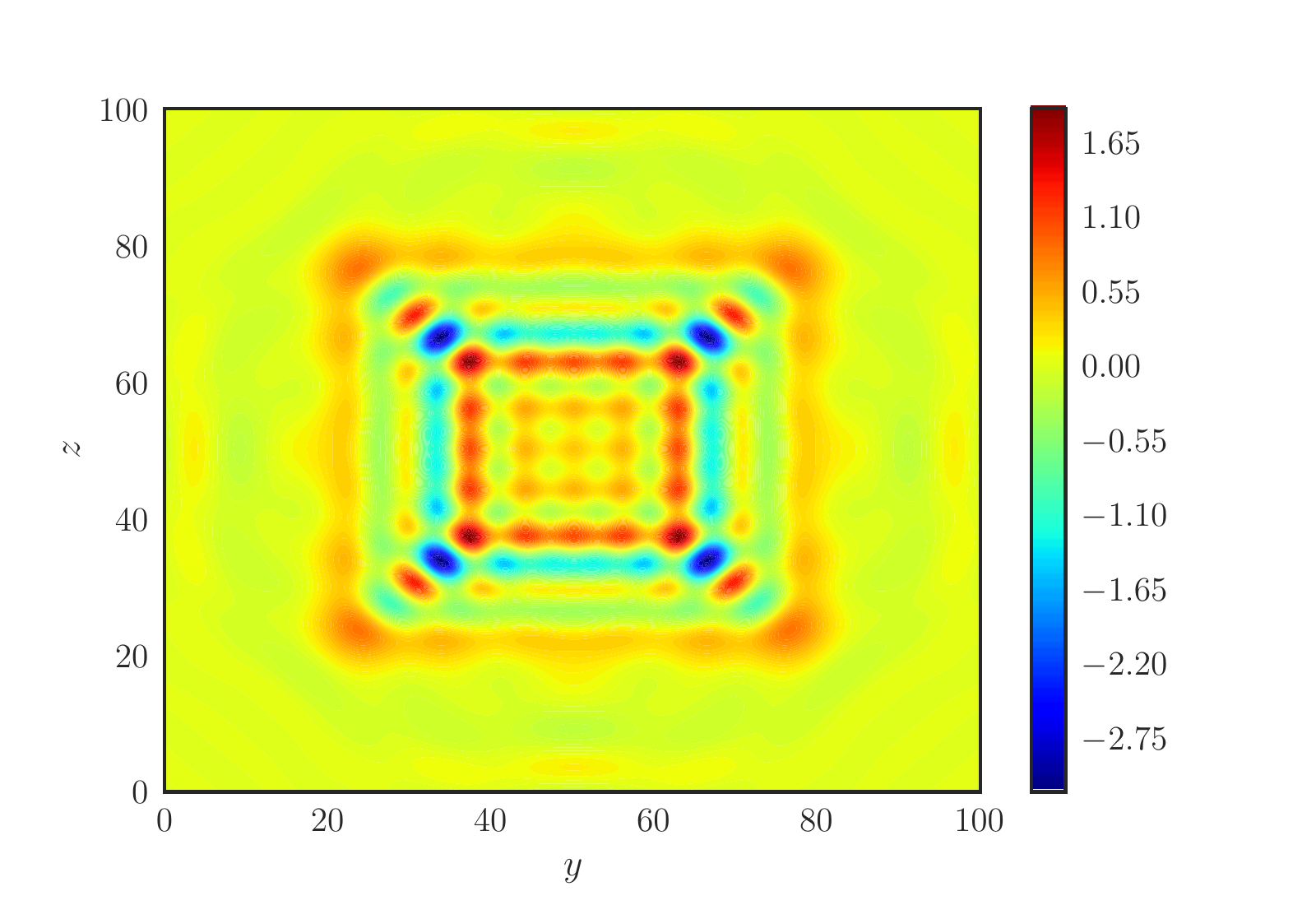}
		     	\caption{Section $x = 17.5 \pi.$}		  		     	
	\end{subfigure}%
	\caption{Scattering on the smoothed cube of size $8\pi$ in $D = [0,32\pi] \times [0,32\pi] \times [0,32\pi]$ for $k=1$. Effective rank of the solution is equal to $81.$ }
	\label{fig:cube}
\end{figure}
\paragraph{Scattering problems on quasiperiodic structures.}
As a final test we took $q(\textbf{r})$ to be a periodic grid of smoothed cubes (structures like this are extremely suitable for QTT computations). Namely, denoting $q(\textbf{r})$ defined by formula \eqref{eq:smoothed-cube} by $q_{\textbf{r}',a}(\textbf{r})$ we solve the equation \eqref{eq:LS} in the domain $[0,2] \times [0,1] \times [0,10]$ on the grid $64 \times 64 \times 1024$ for $q(\textbf{r})$ defined as
$$q(\textbf{r}) = \sum_{i=1}^{19} \sum_{j=0}^3 \sum_{k=0}^1 q_{(0.25+0.5j,0.25+0.5k,0.5i),0.1}(\textbf{r})$$ and for $k = 4 \pi$. $\phi^{\mathrm{inc}}$ in this case is a plane wave propagating in $z$-direction $$\phi^{\mathrm{inc}}(x,y,z) = e^{4\pi\mi z}.$$ Figure~\ref{fig:crystal} demonstrates our results.
\begin{figure}[h!]

		     	\includegraphics[width=0.9\linewidth]{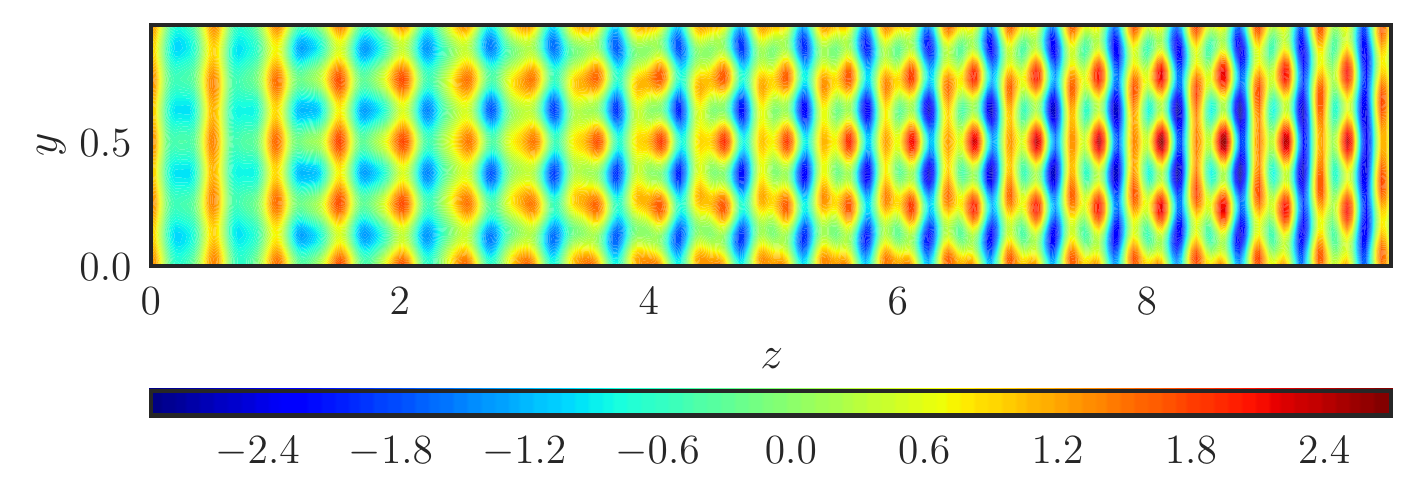}		 		     	

	\caption{Scattering on the periodic cube structure in $D = [0,2] \times [0,1] \times [0,10]$ for $k = 2\pi$. Effective rank of the solution is equal to $67$. Section $x=1$.}
	\label{fig:crystal}
\end{figure}


\section*{Acknowledgements}
This study was supported by the Ministry of Education and Science of the Russian Federation
(grant 14.756.31.0001), by RFBR grants 16-31-60095-mol-a-dk, 16-31-00372-mol-a and by Skoltech NGP program.

\bibliography{tensor,our,vgftt}
\bibliographystyle{ieeetr}

\end{paper}

\end{document}